\documentclass[titlepage,11pt]{article}
\usepackage{latexsym}
\usepackage{amsfonts}
\usepackage{amsmath}
\newcommand\blackslug{\hbox{\hskip 1pt \vrule width 4pt height 8pt depth 1.5pt
        \hskip 1pt}}
\newcommand\bbox{\hfill \quad \blackslug \bigbreak}
\def\d{\hbox{-}}
\def\c{\hbox{-}\cdots\hbox{-}}
\def\l{,\ldots,}

\title{Induced subgraphs of graphs with large chromatic number.\\
II. Three steps towards Gy\'arf\'as' conjectures}
\author{Maria Chudnovsky\thanks{Supported by NSF grants DMS-1001091 and
IIS-1117631.}\\
Princeton University, Princeton, NJ 08544, USA
\\
\\
Alex Scott\\
Oxford University, Oxford, UK
\\
\\
Paul Seymour\thanks{Supported by ONR grant N00014-10-1-0680 and NSF
grant DMS-1265563.}\\
Princeton University, Princeton, NJ 08544}

\date{September 10, 2014; revised \today}

\newtheorem{thm}{}[section]

\newcommand{\Proof}{\noindent{\bf Proof.}\ \ }

\begin{document}
\maketitle
\begin{abstract}
Gy\'arf\'as conjectured in 1985 that for all $k,\ell$, every graph with no clique of size more than $k$ 
and no odd hole of length more than $\ell$ has chromatic number
bounded by a function of $k,\ell$. We prove three weaker statements:
\begin{itemize}
\item Every triangle-free graph with
sufficiently large chromatic number has an odd hole of length different from five;
\item For all $\ell$, every triangle-free graph with sufficiently large chromatic number contains either a $5$-hole
or an odd hole of length more than $\ell$;
\item For all $k,\ell$, every graph with no clique of size more than $k$ and sufficiently large chromatic number 
contains either a $5$-hole
or a hole of length more than $\ell$.
\end{itemize}
\end{abstract}

\section{Introduction}

All graphs in this paper are finite, and without loops or parallel edges. A {\em hole} in a graph $G$ is an induced subgraph
which is a cycle of length at least four, and an {\em odd hole} means a hole of odd length. (The {\em length} of a path
or cycle is the number of edges in it, and we sometimes call a hole of length $n$ an {\em $n$-hole}.) 
In 1985, A. Gy\'arf\'as~\cite{gyarfas} made a sequence of three famous conjectures:

\begin{thm}\label{gyarfasconj0}
{\bf Conjecture:} For every integer $k$ there exists $n(k)$ such that every graph $G$ with no clique of cardinality
more than $k$ and no odd hole has chromatic number at most $n(k)$.
\end{thm}

\begin{thm}\label{gyarfasconj1}
{\bf Conjecture:} For all integers $k,\ell$ there exists $n(k,\ell)$ such that every graph $G$ with no clique of cardinality
more than $k$ and no hole of length more than $\ell$ has chromatic number at most $n(k,\ell)$.
\end{thm}

\begin{thm}\label{gyarfasconj}
{\bf Conjecture:} For all integers $k,\ell$ there exists $n(k,\ell)$ such that every graph $G$ with no clique of cardinality
more than $k$ and no odd hole of length more than $\ell$ has chromatic number at most $n(k,\ell)$.
\end{thm}

In a recent paper~\cite{oddholes}, two of us proved the first conjecture.   
Note that the first two conjectures are special cases of the third.  In the case of the third conjecture,
we might as well assume that $k\ge 2$, and $\ell\ge 3$ and is odd.
Thus it follows from~\cite{oddholes} that conjecture~\ref{gyarfasconj} 
holds for all pairs $(k,\ell)$ when $\ell=3$.
No other cases
have been settled at the time of writing this paper, and the cases when $k=2$ are presumably the simplest to attack next. 
Here we settle the first open case,
when $k=2$ and $\ell=5$. (Since this paper was submitted for publication, we have
proved the second conjecture~\cite{longholes}, and two of us proved the third~\cite{holeseq} 
when $k=2$; part of the proof of the latter uses results of this paper, however,
so this paper is not completely redundant.) 

The conjecture \ref{gyarfasconj} when $(k,\ell) = (2,5)$ asserts that all pentagonal graphs have bounded chromatic number, where we say 
a graph is {\em pentagonal} if every induced odd cycle in it has length five (and in particular, it has no
triangles). Pentagonal graphs might all be $4$-colourable
as far as we know (the $11$-vertex Gr\"{o}tzsch graph is pentagonal and not $3$-colourable), but at least 
they do indeed all have bounded chromatic number.
The following is our main result:

\begin{thm}\label{mainthm}
Every pentagonal graph is $58000$-colourable.
\end{thm}

The proof of \ref{mainthm} occupies almost the whole paper. (Much of the proof needs just that 
$G$ is triangle-free and
has no odd hole of length more than $\ell$, for any fixed $\ell$, and so we have written it in this generality wherever we could.)
We prove:
\begin{itemize}
\item if $G$ has no triangle and no odd hole of length more than $\ell$, and
for every vertex
$v$ the set of vertices with distance at most two from $v$ has chromatic number at most some $k$, then $\chi(G)$ is 
bounded by a function of $k$ and $\ell$; 
\item if $G$ is pentagonal, and $\chi(G)$ is large, then there is an induced subgraph
with large chromatic number in which for every vertex $v$ the set of vertices with distance at most two from $v$ has 
chromatic number at most $5$.
\end{itemize}
Together these imply that every pentagonal graph has bounded chromatic number.
Both of these are consequences of a lemma, a variant of a theorem of~\cite{oddholes}, asserting roughly that for all $\ell$, 
if $G$ is triangle-free
and has no odd hole of length more than $\ell$, and $\chi(G)$ is large, then there is an induced subgraph $H$ such that
for some vertex $v_0$ of $H$, if we partition $V(H)$ by distance in $H$ from $v_0$, then all these ``level sets'' are stable except for
one with large $\chi$. We prove this lemma first, and then apply it to prove the two bulleted statements in later sections.

At the end of this paper, we prove two further special cases of conjectures \ref{gyarfasconj1} and \ref{gyarfasconj}: we show
that conjecture \ref{gyarfasconj1} holds if in addition we assume that $G$ contains no $5$-hole, and that
conjecture \ref{gyarfasconj} holds if in addition we assume that $G$ contains no triangle and no $5$-hole
More precisely, we prove
the next two results, where $\omega(G)$ denotes the size of the largest clique of $G$:

\begin{thm}\label{longoddhole}
Let $\ell\ge 2$ be an integer, and let $G$ be a triangle-free graph with no $5$-hole and no odd hole of length more than $2\ell+1$.
Then $\chi(G)\le (\ell+1)4^{\ell-1}$.
\end{thm}

\begin{thm}\label{longhole}
Let $\ell\ge 3$ be an integer, and let $G$ be a graph with no $5$-hole and no hole of length more than $\ell$. Then
$$\chi(G)\le (2\ell-2)^{2^{\omega(G)}}.$$
\end{thm}
The last was proved (but not published) by the second author some time ago, and improves on~\cite{scott}.

\section{Lollipops}

In \cite{gyarfas}, Gy\'arf\'as gave a neat proof that for any fixed path $P$, all graphs with no induced subgraph isomorphic to $P$ and with bounded clique number
also have bounded chromatic number, and in this section we use basically the same proof for a lemma that we need later.
If $X\subseteq V(G)$, the subgraph of $G$ induced on $X$ is denoted by $G[X]$, and we sometimes write
$\chi(X)$ for $\chi(G[X])$ when there is no danger of ambiguity. 
If $x\in V(G)$ and $Y\subseteq V(G)$, the {\em distance} in $G$ of $x$ from $Y$ (or of $Y$ from $x$) is the length of the shortest path 
containing $x$ and a vertex in $Y$.
Let us say a {\em lollipop} in a graph $G$ is a pair $(C,T)$ where 
$C\subseteq V(G)$ and $T$ is an induced path of $G$ with vertices $t_1\c t_k$ in order, say, with $k\ge 2$, satisfying:
\begin{itemize}
\item $V(T)\cap C=\emptyset$;
\item $G[C]$ is connected;
\item  $t_k$ has a neighbour in $C$; and
\item $t_1\l t_{k-1}$ have no neighbours in $C$.
\end{itemize}
With this notation, the {\em cleanliness} of a lollipop $(C,T)$ in $G$ is the maximum $l$ such that $t_1\l t_l$ all have distance (in $G$) 
at least three from $C$ 
(or $0$ if $t_1$ has distance two from $C$). It follows that the cleanliness is at most $k-2$.
We call $t_1$ the {\em end} of the lollipop.
If $(C,T)$ and $(C',T')$ are lollipops in $G$, we say the second is a {\em licking} of the first if $C'\subseteq C$, and
they have the same end, and $T$ is a subpath of $T'$, and $V(T')\subseteq V(T)\cup C$ 
(and consequently the cleanliness of $(C',T')$ is at least that of
$(C,T)$). We observe first:

\begin{thm}\label{licking} Let $(C,T)$ be a lollipop in $G$, and let $C'\subseteq C$ be non-null, such that $G[C']$ is connected. Then
there is a path $T'$ of $G$ such that $(C',T')$ is a licking of $(C,T)$.
\end{thm}
\Proof
Let $T$ be $t_1\c t_k$, where $(C,T)$ has end $t_1$. If $t_k$ has a neighbour in $C'$ then we may take $T=T'$, so we assume
not. Since $t_k$ has a neighbour in $C$, there is a path $P$ of $G$ 
with one end $t_k$, and with $V(P)\subseteq C\cup \{t_k\}$, such that the other end of $P$ has a neighbour in $C'$.
Choose such a path $P$ with minimum length. Then $V(P)\cap C'=\emptyset$, and $P'=T\cup P$ is an induced path. No vertex of $P'$ 
has a neighbour in $C'$ except its last, and
so $(C', T\cup P)$ is a licking of $(C,T)$ as required.
This proves \ref{licking}.~\bbox

For a vertex $v$ of $G$, we denote the set of neighbours of $v$ in $G$ by $N(v)$ or $N_G(v)$, and for $r\ge 1$, we denote the set
of vertices at distance exactly $r$ 
from $v$ by $N^r(v)$ or $N^r_G(v)$. We need the following:

\begin{thm}\label{lollipop}
Let $h,\kappa\ge 0$ be integers. Let $G$ be a graph such that $\chi(N^2(v))\le \kappa$ for every vertex $v$; and let
$(C,T)$ be a lollipop in $G$, with $\chi(C)> h\kappa$.
Then there is a licking $(C',T')$ of $(C,T)$, with cleanliness at least $h$ more than the cleanliness of $(C,T)$, and such that
$\chi(C')\ge \chi(C)-h\kappa$.
\end{thm}
\Proof We proceed by induction on $h$. If $h=0$ we may take $(C',T')=(C,T)$; so we assume that $h>0$, and that the result holds for 
$h-1$. Let $(C,T)$ have cleanliness $c$ say (where possibly $c=0$), 
and let $T$ have vertices $t_1\c t_k$ in order, where $t_1$ is the end. 
Thus $t_i$ has distance at least three from $C$ for $1\le i\le c$, and so $k\ge c+2$. Since $\chi(N^2(t_{c+1}))\le \kappa$,
it follows that 
$\chi(C\setminus N^2(t_{c+1}))\ge \chi(C)-\kappa$,
and so there is a component $C''$ of $C\setminus N^2(t_{c+1})$ with 
$$\chi(C'')\ge \chi(C)-\kappa>(h-1)\kappa\ge 0.$$
By \ref{licking}, there exists $T''$ such that $(C'',T'')$ is a licking of $(C,T)$. 
Since $t_{c+1}$ has distance at least three from $C''$, it follows
that $(C'',T'')$ has cleanliness at least $c+1$. From the inductive hypothesis, there is a licking $(C',T')$ of $(C'',T'')$
and hence of $(C,T)$ that satisfies the theorem. This proves \ref{lollipop}.~\bbox

\section{Stable levelling}

Let $G$ be a graph. A {\em levelling} $\mathcal{L}$ in $G$ is a sequence $(L_0, L_1\l L_k)$ of disjoint subsets of $V(G)$,
with the following properties:
\begin{itemize}
\item $|L_0|=1$;
\item for each $i$ with $1\le i\le k$, every vertex in $L_i$ has a neighbour in $L_{i-1}$; and
\item for $0\le i,j\le k$ with $|j-i|>1$, there are no edges between $L_i$ and $L_j$.
\end{itemize}

The levelling $\mathcal{L}$ is called {\em stable} if each of the sets
$L_0\l L_{k-1}$ is stable (we do not require $L_k$ to be stable).
For $1\le i\le k$, a {\em parent} of $v\in L_i$ is a neighbour $u$
of $v$ in $L_{i-1}$ (and we also say $v$ is a {\em child} of $u$). 

The next result is a variant of a theorem
proved in~\cite{oddholes}; we could use that theorem directly, but the modification here works better numerically.
Let the {\em odd hole number} of $G$ be the length of the longest induced odd cycle in $G$ (or 1, if $G$ is bipartite).
If $(L_0\l L_k)$ is a stable levelling, we call $L_k$ its {\em base}. 

\begin{thm}\label{parentrule}
Let $G$ be a triangle-free graph with odd hole number at most $2\ell+1$, such that $\chi(N^2(v))\le \kappa$ for every vertex $v$.
Let $(L_0,L_1\l L_k)$ be a levelling in $G$. 
Then there is a stable levelling in $G$ with base of chromatic number at least
$(\chi(L_k) -(\ell-1)\kappa)/2$.
\end{thm}

\Proof We may assume $\ell\ge 1$, since otherwise $G$ is bipartite and the result is trivial. 
Also we may assume that $\chi(L_k)> (2\ell-1)\kappa$, because otherwise the stable levelling $(L_0,L_1)$ satisfies the theorem.
We proceed by induction on $|V(G)|$, and so we may assume:
\begin{itemize}
\item $V(G)=L_0\cup L_1\cup \cdots\cup L_k$;
\item $G[L_k]$ is connected; and 
\item for $0\le i<k$ and every vertex $u\in L_i$, there exists $v\in L_{i+1}$ such that $u$ is its only parent (for if not, we may replace
$L_i$ by $L_i\setminus \{u\}$).
\end{itemize}
Let $L_0=\{s_0\}$, and inductively for $1\le i\le k$, choose $s_i\in L_i$ such that $s_{i-1}$ is its only parent. Then $s_0\d s_1\c s_k$
is an induced path $S$ say.

Now $s_{k-2}$ has no neighbour in $L_k$, so $(L_k,s_{k-2}\d s_{k-1})$ is a lollipop. By \ref{lollipop}, there is a 
licking of this lollipop, say $(C',T')$, with cleanliness at least $2\ell-1$ and with $\chi(C')\ge \chi(L_k)-(2\ell-1)\kappa$.
Let the first $2\ell-1$ vertices of $T'$ be $s_{k-2}\d s_{k-1}\d  t_1\c t_{2\ell-3}$.

Let $N(S)$ be the set of vertices of $G$ not in $S$ but with a neighbour in $S$.
If $v\in L_i\cap N(S)$, then $v$ is adjacent to exactly one of $s_i,s_{i-1}$ and has no other neighbour in $S$; because
every neighbour of $v$ belongs to one of $L_{i-1},L_i,L_{i+1}$, and $G$ is triangle-free, and
$v$ is not adjacent to $s_{i+1}$ since $s_i$ is the only parent of $s_{i+1}$.
So every vertex in $L_i\cap N(S)$ has one of two possible types. 
We say the {\em type} of a vertex $v\in L_i\cap N(S)$ is $\alpha$ where
$\alpha=1$ or $2$ depending whether $v$ is adjacent to $s_{i-1}$ and not to $s_i$, or
adjacent to $s_i$ and not to $s_{i-1}$.

Let us fix a type $\alpha$. Let $V(\alpha)$ be the minimal subset of $V(G)\setminus V(S)$ such that
\begin{itemize}
\item every vertex in $N(S)$ of type $\alpha$ belongs to $V(\alpha)$; and
\item for every vertex $v\in V(G)\setminus (V(S)\cup N(S))$, if some parent of $v$ belongs to $V(\alpha)$ then $v\in V(\alpha)$.
\end{itemize}
Consequently, for every vertex $v\in V(\alpha)$, there is a path starting at $v$ and ending at some vertex in $N(S)$
of type $\alpha$, such that each vertex of the path (except $v$) is the parent of the previous vertex, and no vertex of the path
belongs to $N(S)$ except the last.

There are only two types $\alpha$, and so there is a type $\alpha$ such that
$\chi(V(\alpha)\cap C')\ge \chi(C')/2>0$. Let $C$ be the vertex set of a component of
$G[V(\alpha)\cap C']$
with maximum chromatic number, so 
$$\chi(C)\ge \chi(C')/2\ge  (\chi(L_k)-(2\ell-1)\kappa)/2.$$ 
By \ref{licking}, there is a path $T$ such that $(C,T)$ is a licking of $(C',T')$.

Let $J_k= C$, and for $i=k-1,k-2\l 1$ choose $J_{i}\subseteq V(\alpha)\cap L_{i}$ minimal such that every vertex
in $J_{i+1}\setminus N(S)$ has a neighbour in $J_i$. It follows from the cleanliness of $(C',T')$
that $J_{k-1}\cap N(S)=\emptyset$, and
no vertex in $J_{k-1}$ is adjacent to
any of $s_{k-2},s_{k-1},t_1\l t_{2\ell-3}$.
\\
\\
(1) {\em For $1\le i\le k-1$, if $v\in J_i$ and $v$ is nonadjacent to $s_i$, then there is an induced path $P_v$
between $v$ and $s_i$ of length at least 
$2l-3+2(k-i)$
with interior in $L_{i+1}\cup \cdots\cup L_k$, such that
\begin{itemize}
\item if $i\le k-2$, no vertex in $J_i$
different from $v$ has a neighbour in the interior of $P_v$
\item if $i=k-1$, and $u\in J_i\setminus \{v\}$ has a neighbour in the interior of $P_v$, then 
the induced path between
$u,s_{k-1}$ with interior in $V(P_v)$ has length at least $2\ell-1$.
\end{itemize}}
\noindent Since $v\in J_i$, $v$ has a neighbour in $J_{i+1}\setminus N(S)$ with no other parent in $J_i$;
and so there is a path $v=p_i\d p_{i+1}\c p_k$ such that
\begin{itemize}
\item $p_j\in J_j$ for $i\le j\le k$
\item $p_j\notin N(S)$ for $i<j\le k$
\item $p_{j-1}$ is the only parent of $p_j$ in $J_{j-1}$ for $i<j\le k$.
\end{itemize}
Since $p_{k-1}\in J_{k-1}$, and no vertex in $J_{k-1}$ is adjacent to 
any of $s_{k-2},s_{k-1},t_1\l t_{2l-3}$, it follows that there is an induced path from
$p_{k-1}$ to $s_{k-1}$ with interior in $L_k$ containing all of $t_1\l t_{2l-3}$ and at least one more vertex of $L_k$,
and therefore with length at least $2\ell-1$. Its union with the path $p_i\c p_{k-1}$ and the path $s_{k-1}\d s_{k-2}\c s_i$ is an induced path
between $v$ and $s_i$, of length at least $2\ell-3+2(k-i)$. If $u\in J_i\setminus \{v\}$
and has a neighbour in the interior of $P_v$, then since $u$ is nonadjacent to all of $s_{i+1}\l s_{k-1},p_{i+1}\l p_{k-1}$
(because $u$ has no neighbour in $L_{i+2}\cup \cdots\cup L_k$, and $s_{i+1}$ has a unique parent $s_i$, and $p_{i+1}$ has no parent
in $J_i$ except $p_i$), it follows that $i=k-1$; and since no vertex in $J_{k-1}$ is adjacent to   
any of $s_{k-2},s_{k-1},t_1\l t_{2l-3}$, this proves (1).

\bigskip

For $1\le i\le k$ and for every vertex $v\in J_i$, either $v\in N(S)$ or it has a parent in $J_{i-1}$;
and so there is a path $v=r_i\d r_{i-1}\c r_h$
for some $h\le i$, such that $r_j\in J_j$ for $h\le j\le i$, and $r_h\in N(S)$, and $r_j\notin N(S)$ for $h+1\le j\le i$.
Since $r_h$ has a neighbour in $S$, one of
$$r_i\d  r_{i-1}\c r_h\d s_{h-1}\d  s_h\d s_{h+1}\c s_i,$$
$$r_i, r_{i-1}\c r_h\d s_h\d s_{h+1}\c s_i$$
is an induced path (the first if $\alpha=1$ and the second if $\alpha=2$). We choose some such path and call it $R_v$.
Note that for all $v\in J_1\cup\cdots\cup J_k$, the path $R_v$ has even length if $\alpha=1$, and odd length otherwise.
\\
\\
(2) {\em For $0\le i\le k-1$, $J_i$ is stable.}
\\
\\
Suppose that $u,v\in J_i$ are adjacent. Since $G$ is triangle-free and $u,v$ have the same type, not both $u,v\in N(S)$.
Suppose that $u\in N(S)$, and hence $v\notin N(S)$. Since $N(S)\cap J_{k-1}=\emptyset$
it follows that $i\le k-2$.
Consequently $u$ has no neighbour in the interior of $P_v$, where $P_v$ is as in (1), and so $P_v\cup R_v$,
$s_i\d P_v\d v\d u\d R_u\d s_i$ are both holes of length at least $2\ell+2$, of different parity, which is impossible.
So $u,v\notin N(S)$.
We claim that there is a path $P$ of length at least $2\ell-1$, from one of $u,v$ to $s_i$, with interior in $L_{i+1}\cup\cdots\cup L_k$,
such that the other (of $u,v$) has no neighbour in its interior. For if $u$ has no neighbour in the interior of $P_v$
then we may take $P=P_v$, where $P_v$ is as in (1); and if $u$ has such a neighbour, let $P$ be the induced path
between $u$ and $s_i$ with interior a subset of the interior of $P_v$. Note that in the second case, $v$ has no neighbour 
in the interior of $P$, since $G$ is triangle-free.
This proves that the desired path $P$ exists; say from $v$ to $s_i$. Now the union of $P$ and $R_v$ is a hole of length at least
$2\ell+2$, and so $P,R_v$ have the same parity. But the union of $P$ and the path $v\d u\d R_u\d s_i$ is also a hole, of length
at least $2\ell+3$, and since $R_u,R_v$ have the same parity this is impossible. This proves (2).

\bigskip
If $\alpha=1$ let $M_i=\{s_i\}\cup J_i$ for $0\le i\le k$, and if 
$\alpha=2$ let $M_0=\{s_1\}, M_i=\{s_{i+1}\}\cup J_i$ for $1\le i<k$, and $M_k= J_k$. In each case $(M_0\l M_k)$
is a levelling satisfying the theorem.
This proves \ref{parentrule}.~\bbox

We deduce:

\begin{thm}\label{getlevelling}
Let $G$ be pentagonal, and let $n\ge 1$ be an integer. If $\chi(G)\ge 10n-9$, there is a stable levelling in $G$ with base of 
chromatic number at least $n$.
\end{thm}
\Proof
Let $G'$ be a component of $G$ with $\chi(G')=\chi(G)$. Choose $v_0\in V(G')$, and for $i\ge 0$ let $L_i$
be the set of vertices in $G'$ with distance $i$ from $v_0$. There exists $k$ such that $\chi(L_k)\ge \chi(G)/2$
and hence $\chi(L_k)\ge 5n-4$.
Now $(L_0\l L_k)$ is a levelling in $G$. By \ref{parentrule}, taking $\ell=2$ and $\kappa= n-1$,
either 
\begin{itemize}
\item there is a vertex $v$ with $\chi(N^2(v))\ge n$, and hence there is a levelling $(M_0, M_1, M_2)$ with $\chi(M_2)\ge n$, 
necessarily stable, or 
\item there is a stable levelling $(M_0\l M_{k})$ in $G$ with 
$\chi(M_k)\ge (\chi(L_k)- 3(n-1))/2\ge n-1/2$. 
\end{itemize}
In either case the theorem holds.~\bbox

\section{Reducing to bounded radius}

Let $(L_0\l L_k)$ be a levelling. If $0\le i\le j\le k$ and $u\in L_i$ and $v\in L_j$, and there
is a path between $u,v$ of length $j-i$ with one vertex in each of $L_i,L_{i+1}\l L_j$, we say that $u$ is an {\em ancestor} of $v$
and $v$ is a {\em descendant} of $u$. 

\begin{thm}\label{boundedrad}
Let $G$ be a triangle-free graph with odd hole number at most $2\ell+1$.
For $r = 2,3$, let $\chi(N^r(v))\le \kappa_r$ for every vertex $v$.
Then $\chi(G) \le (12\ell-6)\kappa_2+ 4\kappa_3+8$.
\end{thm}
\Proof Suppose that $\chi(G) > (12\ell-6)\kappa_2+ 4\kappa_3+8$.
There is a levelling in $G$ with base of chromatic number at least $\chi(G)/2$, and so 
by \ref{parentrule}, there is a stable levelling 
$(L_0\l L_k)$ in $G$ with 
$$\chi(L_k)\ge \chi(G)/4 -(\ell-1/2)\kappa_2> (2\ell-1)\kappa_2+ \kappa_3+2.$$ 
We may choose it in addition such that $G[L_k]$
is connected, and for $0\le i<k$ every vertex in $L_i$ has a descendant in $L_k$. Since
$\chi(L_k)>1$ it follows that $k>1$.
Choose $a_{k-2}\in L_{k-2}$. Let $X_1$ be the set of descendants of $a_{k-2}$ in $L_k$;
thus $X_1\ne \emptyset$, and $\chi(X_1)\le \kappa_2$, and since $\chi(L_k)>\kappa_2$, 
there is a component $C_1$ of $G[L_k\setminus X_1]$ with 
$$\chi(C_1)\ge \chi(L_k)-\kappa_2> (2\ell-2)\kappa_2+ \kappa_3+2.$$
Since $G[L_k]$ is connected and $X_1\ne \emptyset$, there exists $a_k\in X_1$ with a neighbour in $C_1$.
Let $a_{k-1}$ be a parent of $a_k$ and child of $a_{k-2}$. 

Let $X_2$ be the set of neighbours of $a_k$ in $C_1$; then $X_2$ is stable and nonempty, and since
$\chi(C_1)>1$, there is a component
$C_2$ of $C_1\setminus X_2$ with 
$$\chi(C_2)\ge \chi(C_1)-1> (2\ell-2)\kappa_2+ \kappa_3+1,$$
and a neighbour $b_k\in L_k$ of $a_k$ with a neighbour in $C_2$.
Let $b_{k-1}$ be a parent of $b_k$. Thus $b_{k-1}, a_{k-2}$ are nonadjacent since $X_1\cap C_1=\emptyset$. Also $b_{k-1}, a_{k-1}$
are nonadjacent since $L_{k-1}$ is stable, and $b_{k-1},a_k$ are nonadjacent since $G$ is triangle-free, and similarly
$a_{k-1}, b_k$ are nonadjacent. Consequently
$a_{k-2}\d a_{k-1}\d a_k\d b_k\d b_{k-1}$ is an induced path of $G$. 

Let $X_3$ be the set of all children of $b_{k-1}$; then since $X_3$ is stable, and $\chi(C_2)>1$, it follows that there is a component
$C_3$ of $C_2\setminus X_3$ with 
$$\chi(C_3)\ge \chi(C_2)-1> (2\ell-2)\kappa_2+ \kappa_3,$$ 
and a child $c_k$ of $b_{k-1}$ with a neighbour in $C_3$,
taking $c_k=b_k$ if $b_k$ has a neighbour in $C_3$.
Thus $(C_3,b_{k-1}\d c_k)$ is a lollipop. By \ref{lollipop}, since $\chi(C_3)>  (2\ell-2)\kappa_2$,
there is a licking $(C_4,T)$ of $(C_3,b_{k-1}\d c_k)$, with cleanliness at least $2\ell-2$, such that
$$\chi(C_4)\ge \chi(C_3)-(2\ell-2)\kappa_2> \kappa_3.$$ 
Let $T$ have vertices $t_1\d t_2\d t_3\c t_m$ say, where
$m\ge 2\ell$ and $t_1=b_{k-1}$ and $t_2 = c_k$. Note that if $b_k\ne c_k$ then $b_k$ has no neighbour in $C_3$ and in particular
$b_k$ has no neighbour in $T$ except $t_1$.

Let $X_4$ be the set of all vertices of $C_4$ with distance three from $b_{k-1}$.
Since $\chi(X_4)\le \kappa_3$, and $\chi(C_4)-\kappa_3>0$, 
there is a component $C_5$ of $C_4\setminus X_4$.
By \ref{licking}, there is a licking $(C_5,S)$ say of $(C_4,T)$. Let $S$ have vertices 
$t_1\c t_n$ 
say where $n\ge m$. Let $t_{n+1}\in V(C_5)$ be adjacent to $t_n$, and let $d_{k-1}$ be a parent of $t_{n+1}$. 
Choose $i$ with $1\le i\le n+1$
minimum such that $d_{k-1}$ is adjacent to $t_i$. Note that $d_{k-1}$ is nonadjacent to all of $t_1\l t_{2\ell-2}$
since $(C_4,T)$ has cleanliness at least $2\ell-2$ and hence so does $(C_5,S)$; and so $i>2\ell-2$. Let $P$ be the path
$t_2\d t_3\c t_i\d d_{k-1}$. This path $P$ is induced and has length $i-1\ge 2\ell-2$.

Choose parents $b_{k-2}, d_{k-2}$ of $b_{k-1}, d_{k-1}$. Since
Since $t_{n+1}$ is in $C_5$, it follows that $b_{k-1},d_{k-1}$ have distance at least three; and consequently
$b_{k-2}\ne d_{k-2}$, and 
$b_{k-2}$ is nonadjacent to $d_{k-1}$, and $b_{k-1}$ is nonadjacent to $d_{k-2}$.
Now $a_{k-2}\ne d_{k-2}$, since $a_{k-2}$ has no descendant in $C_1$, and 
$d_{k-2}$ has a descendant $t_{n+1}$ in $C_5$ and hence in $C_1$. For the same reason $a_{k-2}$
is nonadjacent to $d_{k-1}$, and in particular $a_{k-1}\ne d_{k-1}$.

Since $L_0\l L_{k-3}$ are stable, there is an induced path
between $b_{k-2}, d_{k-2}$ of even length with interior in $L_0\cup\cdots\cup L_{k-3}$, and its union with
the path $b_{k-2}\d b_{k-1}\d t_2\d P\d d_{k-1}\d d_{k-2}$ is a hole of length at least $2\ell+3$, which consequently has even length; and so 
$P$ has odd length. Now there is an even induced path $Q$ between $a_{k-1}, d_{k-1}$ with interior in $L_0\cup\cdots\cup L_{k-2}$, 
not containing any neighbour of $b_{k-1}$; for if $a_{k-1}, d_{k-2}$ are adjacent then the path $a_{k-1}\d d_{k-2}\d d_{k-1}$
satisfies our requirements, and otherwise any even induced path between $a_{k-2}, d_{k-2}$ with interior in $L_0\cup \cdots\cup L_{k-3}$
(extended by the edges $a_{k-1}a_{k-2}$ and $d_{k-1}d_{k-2}$) provides the desired path. If
$b_k\ne c_k$ then
$$a_{k-1}\d a_k\d b_k\d b_{k-1}\d c_k\d P\d d_{k-1}\d Q\d a_{k-1}$$
is an odd hole of length at least $2\ell+4$, while if $b_k=c_k$ then
$$a_{k-1}\d a_k\d b_k\d P\d d_{k-1}\d Q\d a_{k-1}$$
is an odd hole of length at least $2\ell+2$, in either case a contradiction. This proves \ref{boundedrad}.~\bbox

Next we prove a variant of \ref{boundedrad} in which $\kappa_3$ is eliminated, the following.
The proof is almost the same, but differs in a couple of key places, and we felt it best to write it out completely,
despite the duplication.

\begin{thm}\label{2rad}
Let $G$ be a triangle-free graph with odd hole number at most $2\ell+1$.
Let $\chi(N^2(v))\le \kappa_2$ for every vertex $v$.
Then $\chi(G)\le (40\ell+28)\kappa_2+40$.
\end{thm}
\Proof Suppose that $\chi(G) >  (40\ell+28)\kappa_2+40$.
Choose a levelling $(M_0\l M_k)$ in $G$ with base of chromatic number at least $\chi(G)/2$, and 
let $H = G[M_k]$. 
\\
\\
(1) {\em There exists $v\in M_k$ and a parent $v_{k-1}\in M_{k-1}$ of $v$, and a levelling $(L_0\l L_3)$ in $H$
with $L_0=\{v\}$, such that $L_3$ is disjoint from $N^1_G(v_{k-1})\cup N^2_G(v_{k-1})$, and
$\chi(L_3) >(2\ell+4)\kappa_2+2.$}
\\
\\
Choose $v\in V(H)$ with $\chi(N^3_H(v))$ maximum;
$\chi(N^3_H(v))=\kappa_3$ say. 
By \ref{boundedrad} applied in $H$, 
$$\chi(H)\le (12\ell-6)\kappa_2+ 4\kappa_3+8,$$
and since $\chi(H)\ge\chi(G)/2$, it follows that
$$\kappa_3\ge \chi(G)/8-(3\ell-3/2)\kappa_2 -2> (2\ell+5)\kappa_2 + 3.$$
Let $v_{k-1}\in M_{k-1}$ be a parent of $v$. Now the set of neighbours of $v_{k-1}$
is stable, and $\chi(N^2_G(v))\le \kappa_2$, and so there exists $L_3\subseteq N^3_H(v)$
disjoint from $N^1_G(v_{k-1})\cup N^2_G(v_{k-1})$ with 
$$\chi(L_3)\ge \kappa_3-\kappa_2-1>(2\ell+4)\kappa_2+2.$$
Thus $(\{v\}, N^1_H(v), N^2_H(v), L_3)$ is the desired levelling in $H$. This proves (1).

\bigskip
It follows that $v_{k-1}$ has no neighbour in $L_1$, since $G$ is triangle-free, and has no neighbour in $L_2\cup L_3$, since
$L_3$ is disjoint from $N^1_G(v_{k-1})\cup N^2_G(v_{k-1})$. In addition we may choose $L_0\l L_3$ such that $H[L_3]$
is connected, and every vertex in $L_1$ has a descendant in $L_3$. 
Note that $L_2$ might not be stable, but $L_1$ is stable since $G$ is triangle-free.
Choose $a_{1}\in L_{1}$. Let $X_1$ be the set of descendants of $a_{1}$ in $L_3$;
thus $X_1\ne \emptyset$, and $\chi(X_1)\le \kappa_2$, and since $\chi(L_3)>\kappa_2$, there is a component $C_1$ of $G[L_3\setminus X_1]$ with
$$\chi(C_1)\ge \chi(L_3)-\kappa_2> (2\ell+3)\kappa_2+2.$$
Since $G[L_3]$ is connected and $X_1\ne \emptyset$, there exists $a_3\in X_1$ with a neighbour in $C_1$.
Let $a_{2}\in L_2$ be adjacent to $a_1,a_3$.

Let $X_2$ be the set of neighbours of $a_3$ in $C_1$; then $X_2$ is stable and nonempty, and since
$\chi(C_1)>1$, there is a component
$C_2$ of $C_1\setminus X_2$ with
$\chi(C_2)\ge \chi(C_1)-1> (2\ell+3)\kappa_2+1,$
and a neighbour $b_3\in L_3$ of $a_3$ with a neighbour in $C_2$.
Let $b_{2}\in L_2$ be adjacent to $b_3$. Thus $b_{2}, a_{1}$ are nonadjacent since $X_1\cap C_1=\emptyset$. 
Also $b_{2},a_3$ are nonadjacent since $H$ is triangle-free. (But 
$b_{2}, a_{2}$ might be adjacent.) Consequently one of
$a_{1}\d a_{2}\d a_3\d b_3\d b_{2}$, $a_{1}\d a_{2}\d b_{2}$ is an induced path of $H$ with even length.

Let $X_3$ be the set of all children of $b_{2}$; then since $X_3$ is stable, and $\chi(C_2)>1$, it follows that there is a component
$C_3$ of $C_2\setminus X_3$ with
$$\chi(C_3)\ge \chi(C_2)-1>(2\ell+3)\kappa_2,$$
and a child $c_3$ of $b_{2}$ with a neighbour in $C_3$,
taking $c_3=b_3$ if $b_3$ has a neighbour in $C_3$.
Thus $(C_3,b_{2}\d c_3)$ is a lollipop. By \ref{lollipop} applied in $G$ (not just in $H$), since $\chi(C_3)>  (2\ell-2)\kappa_2$,
there is a licking $(C_4,T)$ of $(C_3,b_{2}\d c_3)$, with cleanliness at least $2\ell-2$ in $G$, such that
$$\chi(C_4)\ge \chi(C_3)-(2\ell-2)\kappa_2>5\kappa_2.$$
Let $T$ have vertices $t_1\d t_2\d t_3\c t_m$ say, where
$m\ge 2\ell$ and $t_1=b_{2}$ and $t_2 = c_3$. Note that if $b_3\ne c_3$ then $b_3$ has no neighbour in $C_3$ and in particular
$b_3$ has no neighbour in $T$ except $t_1$.

Let $b_1\in L_1$ be adjacent to $b_2$; and
let $X_4$ be the set of all vertices of $C_4$ with distance two in $G$ from one of $a_1,a_2,a_3,b_1,b_3$.
Since $\chi(X_4)\le 5\kappa_2$, and $\chi(C_4)-5\kappa_2>0$,
there is a component $C_5$ of $C_4\setminus X_4$.
By \ref{licking}, there is a licking $(C_5,S)$ say of $(C_4,T)$. Let $S$ have vertices
$t_1\c t_n$
say where $n\ge m$. Let $t_{n+1}\in V(C_5)$ be adjacent to $t_n$, and let $d_{k-1}\in M_{k-1}$ be a parent of $t_{n+1}$.
(Note that $d_{k-1}$ belongs to $M_{k-1}$, not to $L_2$; this is where this proof differs essentially
from the proof of \ref{boundedrad}.)
Choose $i$ with $1\le i\le n+1$
minimum such that $d_{k-1}$ is adjacent to $t_i$. Note that $d_{k-1}$ is nonadjacent to all of $t_1\l t_{2\ell-2}$
since $(C_4,T)$ has cleanliness in $G$ at least $2\ell-2$ and hence so does $(C_5,S)$; and so $i>2\ell-2$. Note also 
that $d_{k-1}$ is nonadjacent to all of $a_1,a_2,a_3,b_1,b_2,b_3$, from the definition of $C_5$.

The path $c_3=t_2\d t_3\c t_i\d d_{k-1}$ is induced and has length $i-1\ge 2\ell-2$. Now 
there is an induced path between $d_{k-1}, v$ with interior a subset of $M_0\cup M_1\cup \cdots\cup M_{k-2}\cup \{v_{k-1}\}$;
and the union of this path with the previous one is an induced path $P$ of length at least $2\ell-1$ between $t_2$ and $v$.
Note that none of $a_1,a_2,a_3,b_1,b_2,b_3$ have neighbours in the interior of $P$. Now the union of $P$ and the path
$c_3\d b_2\d b_1\d v$ is a hole of length at least $2\ell+2$, and so is even; and hence $P$ has odd length.
Let $Q$ be the path
\begin{itemize}
\item $c_3\d b_2\d a_2\d a_1\d v$ if $b_2,a_2$ are adjacent;
\item $c_3\d b_2\d b_3\d a_3\d a_2\d a_1\d v$ if $b_2,a_2$ are nonadjacent and $c_3\ne b_3$; and
\item $c_3\d a_3\d a_2\d a_1\d v$ if $b_2,a_2$ are nonadjacent and $c_3=b_3$.
\end{itemize}
In each case $Q$ is between $c_3,v$, and has even length, at least four.
The union of $P$ and $Q$
is therefore an odd hole of length at least $2\ell+3$, a contradiction. This proves \ref{2rad}.~\bbox

\section{The Gr\"{o}tzsch graph}

Let $G$ be a graph, and $H$ an induced subgraph of $G$. We say a levelling $(L_0\l L_k)$ in $G$ is {\em over $H$} if
$V(H)\subseteq L_k$. For $n\ge 1$ an {\em $n$-covering (in $G$, over $H$)} 
is a sequence of graphs $H=G_0,G_1\l G_n=G$, such that for $1\le i\le n$ there is a stable levelling 
in $G_{i}$ over $G_{i-1}$.
For $n\ge 1$, let us say a graph $H$ is {\em $n$-coverable} if there is 
an $n$-covering over $H$ in some pentagonal graph $G$ (and in particular, $H$ itself is pentagonal).

The {\em Gr\"{o}tzsch graph} has vertex set $\{a_1\l a_5, b_1\l b_5, c\}$, where $a_1\d a_2\c a_5\d a_1$ is a cycle,
$a_i,b_i$ are both adjacent to $a_{i-1}$ and $a_{i+1}$ for $1\le i\le 5$ (reading subscripts modulo $5$), and $c$ is adjacent to 
$b_1\l b_5.$ We call the $5$-hole $a_1\d a_2\c a_5\d a_1$ its {\em rim} and $c$ its {\em apex}. 

\begin{thm}\label{myccover}
The Gr\"{o}tzsch graph is not $1$-coverable.
\end{thm}
\Proof Suppose it is, and let $G$ be pentagonal, with a stable levelling $(L_0\l L_k)$, such that $G[L_k]$ has an induced subgraph
$H$ isomorphic to the Gr\"{o}tzsch graph. Let $V(H)$ be labelled as above. We may assume that $L_k=V(H)$, and $L_{k-1}$
is minimal such that every vertex in $V(H)$ has a neighbour in $L_{k-1}$. For each $v\in L_{k-1}$, let $H(v)$ denote the set of
neighbours of $v$ in $V(H)$. Consequently:
\\
\\
(1) {\em For each $v\in L_{k-1}$, there exists $u\in H(v)$ with no neighbour in $L_{k-1}$ except $v$.}
\\
\\
We call such a vertex $u$ a {\em dependent} of $v$.
If $u,v\in L_{k-1}$, by a $u\d v$ {\em gap} we mean an induced path $P$ of $G$, with one end in $H(u)$ and the other in $H(v)$,
and with no other vertex in $H(u)\cup H(v)$ (a vertex in $H(u)\cap H(v)$ forms a $1$-vertex gap.) Thus a $u\d v$ gap is the interior
of an induced path between $u$ and $v$. 
\\
\\
(2) {\em For all $u,v\in L_{k-1}$, no $u\d v$ gap has length three.}
\\
\\
For suppose some $u\d v$ gap has length three; then there is an induced path between $u,v$ of length five,
with interior in $L_k$. But $u,v$ have neighbours in $L_{k-2}$, and so are joined by an induced path of even length with interior
in the top of the levelling; and the union of these two paths is an odd hole of length at least seven, which is impossible.
\\
\\
(3) {\em For every four-vertex induced path $u_1\d u_2\d u_3\d u_4$ of $H$, if $v,v'\in L_{k-1}$ and $u_1\in H(v)$ and $u_4\in H(v')$,
then either one of $u_1,u_2\in H(v')$, or one of $u_3, u_4\in H(v)$.}
\\
\\
Because $H(v), H(v')$ are stable sets since $G$ is triangle-free; and from (2) this path is not a $u\d v$ gap; and the claim follows.
\\
\\
(4) {\em $|H(v_0)|\ge 2$ for all $v_0\in L_{k-1}$ with $c\in H(v_0)$.}
\\
\\
For suppose that $H(v_0)=\{c\}$. Then by (1), $c$ has no other neighbour in $L_{k-1}$. So for every four-vertex induced path of $H$
ending at $c$, say $u_1\d u_2\d u_3\d c$, and for all $v\in L_{k-1}$ with $u_1\in H(v)$, (3) implies that $u_3\in H(v)$ (because
$u_1,u_2\notin H(v_0)$ since $|H(v_0)|=1$, and $c\notin H(v)$ since $c$ is a dependent of $v_0$). Choose $v_1\in L_{k-1}$
with $a_1\in H(v_1)$. From
$a_1\d a_5\d b_1\d c$ it follows that $b_1\in H(v_1)$, and similarly $b_3,b_4\in H(v_1)$. Since $H(v_1)$ is stable, and the set
$\{a_1,b_1,b_3,b_4\}$ is a maximal stable set of $H$, it follows that $H(v_1)=\{a_1,b_1,b_3,b_4\}$. Choose $v_3\in L_{k-1}$
with $a_3\in H(v_3)$; then (from the symmetry of $H$ taking $a_3$ to $a_1$) it follows that $H(v_3)=\{a_3,b_3, b_1,b_5\}$. 
But $b_4\d a_5\d a_4\d b_5$ is a $v_1\d v_3$
gap contradicting (2). This proves (1).
\\
\\
(5) {\em $|H(v_0)|= 3$ for all $v_0\in L_{k-1}$ with $c\in H(v_0)$.}
\\
\\
No stable set of $H$ containing $c$ has cardinality more than three, so we just need to show that $|H(v_0)|\ne 2$. Suppose
not; then from the symmetry of $H$, we may assume that $H(v_0)=\{c,a_1\}$.
One of $c,a_1$ is a dependent of $v_0$. 

Suppose first that $c$ is a dependent of $v_0$. Choose $v_5\in L_{k-1}$ with
$a_5\in H(v_5)$. From $a_5\d a_4\d b_3\d c$ and (3) it follows that $b_3\in H(v_5)$, and from $a_5\d a_4\d b_5\d c$
that $b_5\in H(v_5)$. Since $a_2,b_3$ are adjacent it follows that $a_2\notin H(v_5)$; choose $v_2\in L_{k-1}$
with $a_2\in H(v_2)$. Then from the symmetry of $H$ exchanging $a_2,a_5$ and fixing $a_1$,
it follows $b_2,b_4\in H(v_2)$. From $a_5\d b_1\d c\d b_2$ and (3) it follows that $b_2\in H(v_5)$ (since
$a_5,b_1\notin H(v_2)$ because they both have neighbours in $H(v_2)$, and $c\notin H(b_5)$ because it is a dependent of $v_1$).
From the same symmetry, $b_5\in H(v_2)$; and so $a_3\notin H(v_2)$ and $a_4\notin H(v_5)$. 
But then $a_5\d a_4\d a_3\d a_2$ is a $v_5\d v_2$ gap, contrary to (2).

This shows that $c$ is not a dependent of $v_0$, and so $a_1$ is its dependent. Choose $v_3\in L_{k-1}$ with $a_3\in H(v_3)$;
then $b_5\in H(v_3)$ from $a_3\d a_4\d b_5\d a_1$, and $a_5\in H(v_3)$ from $a_3\d b_4\d a_5\d a_1$. Now $a_4\notin H(v_3)$; choose
$v_4\in L_{k-1}$ with $a_4\in H(v_4)$, and then similarly $a_2,b_2\in H(v_4)$. But then $a_5\d b_1\d c\d b_2$ is a $v_3\d v_2$ gap,
a contradiction.
This proves (5).

\bigskip

In view of (5) and the symmetry we may assume henceforth that $H(v_0) = \{a_5,a_2,c\}$. One of $a_5,a_2,c$ is a dependent of $v_0$.
Suppose first that $c$ is a dependent of $v_0$. Choose $v_3\in L_{k-1}$ with $a_3\in H(v_3)$; then $b_5\in H(v_3)$ from
$a_3\d a_4\d b_5\d c$, and $b_3\in H(v_3)$ from $a_3\d a_4\d b_3\d c$. Similarly, let $a_4\in H(v_4)$; then $b_2,b_4\in H(v_4)$.
From $b_4\d a_5\d a_1\d b_5$ it follows that $a_5\in H(v_3)$, and similarly $a_2\in H(v_4)$. But then $a_5\d b_1\d c\d b_2$
is a $v_3\d v_4$ gap, a contradiction.

From the symmetry between $a_2, a_5$, we may therefore assume that $a_5$ is a dependent of $v_0$.
Let $b_2\in H(v_2)$; then $a_4\in H(v_2)$ from $b_2\d a_3\d a_4\d a_5$, and $b_4\in H(v_2)$ from $b_2\d a_3\d b_4\d a_5$.
Also, $a_2\in H(v_2)$ from $a_4\d b_5\d a_1\d a_2$. Let $a_3\in H(v_3)$; then $a_1\in H(v_3)$ from $a_3\d b_2\d a_1\d a_5$,
and $c\in H(v_3)$ from $a_1\d b_5\d c\d b_4$. But then $a_4\d a_5\d b_1\d c$ is a $v_2\d v_3$ gap, a contradiction.
This proves \ref{myccover}.~\bbox

\section{Radius two}

In this section we prove a bound on $\chi(N^2(v))$ for $2$-coverable graphs,
to allow us to apply \ref{boundedrad}. We begin with:

\begin{thm}\label{pentcover}
Let $(L_0\l L_k)$ be a stable levelling in a pentagonal graph $G$, and let $P$ be a $5$-hole of $G[L_k]$.
Choose $S\subseteq L_{k-1}$ minimal such that every vertex in $P$ has a neighbour in $S$. Then 
\begin{itemize}
\item $|S| = 3$;
\item we can label the vertices of $P$ as $p_1\c p_5\d p_1$ in order, and label the elements of $S$ as $a,b,c$, such that 
the edges of $G$ between $S$ and $V(P)$ are $ap_1,ap_3,bp_2,bp_4,cp_5$ and possibly $cp_3$;
\item there exists $z\in L_{k-2}$ adjacent to every vertex in $S$.
\end{itemize}
\end{thm}
\Proof
We begin by proving the first two assertions.
Each vertex in $S$ has at most two neighbours in $P$, because its neighbours form a stable set.
Suppose that every vertex in $S$ has exactly two neighbours in $P$. We may assume that $a\in S$ is adjacent to $p_1,p_3$; then
choose $b\in S$ adjacent to $p_2$. It follows that $b$ is adjacent to one of $p_4,p_5$, say $p_4$. Choose $c\in S$ adjacent to $p_5$;
then $c$ might also be adjacent to one of $p_2,p_3$, and from the symmetry we may assume it is not adjacent to $p_2$; and so
$S=\{a,b,c\}$, and the first two assertions of the theorem hold. We may therefore assume that some vertex in $S$, say $c$,
has only one neighbour in $P$, say $p_5$. From the minimality of $S$, no other vertex in $S$ is adjacent to $p_5$. Choose $a\in S$
adjacent to $p_3$. If $a$ has no more neighbours in $P$, then the path $a\d p_3\d  p_2\d p_1\d p_5\d c$ can be completed via an
even path joining $a,c$ with interior in $L_0\cup \cdots\cup L_{k-2}$ to an odd hole of length at least seven, which is impossible.
So $a$ has another neighbour in $P$, and since $a$ is not adjacent to $p_5$ it is adjacent to $p_1$. Similarly, choose $b\in S$
adjacent to $p_2$; then $b$ is also adjacent to $p_4$. From the minimality of $S$, $S=\{a,b,c\}$ and again the 
first two assertions hold.

For the third assertion, choose $Z\subseteq L_{k-2}$ minimal containing a neighbour of each member of $S$.
Suppose that there are distinct $z_1,z_2\in Z$. From the minimality of $Z$, there exist $s_1,s_2\in S$ such that
for $1\le i,j\le 2$, $z_i$ is adjacent to $s_j$ if and only if $i=j$. But from the second assertion of the theorem,
there is a three-edge path joining $s_1,s_2$ with interior in $V(P)$, say $s_1\d p_1\d p_2\d s_2$, where $p_1p_2$
is an edge of $P$. Then $z_1\d s_1\d p_1\d p_2\d s_2\d z_2$ is an induced path, and can be completed to an 
odd hole of length at least seven via an even induced path joining $z_1,z_2$ with interior in $L_0\cup\cdots\cup L_{k-3}$, which
is impossible. Thus $|Z|=1$, and so the third assertion holds.
This proves \ref{pentcover}.~\bbox

We also need the following lemma.
\begin{thm}\label{nbrz}
Let $G$ be pentagonal, and let $(L_0\l L_k)$ be a stable covering in $G$ of a graph $H$.
Let $z\in V(H)$, and $A$ be the set of all vertices $v\in N^2_H(z)$ such that every neighbour of $v$ in $L_{k-1}$
is adjacent to $z$. Then $\chi(A)\le 2$.
\end{thm}
\Proof
Suppose that $\chi(A)>2$; then there is a $5$-hole $P$ of $G[A]$. Choose a minimal subset $S$ of $N^1_H(z)$
such that every vertex in $p$ has a neighbour in $S$; then by \ref{pentcover}
we may assume that
$S=\{a,b,c\}$, where the edges between $S$ and $V(P)$ are 
$ap_1,ap_3,bp_2,bp_4,cp_5$ and possibly $cp_3$. Choose $v\in L_{k-1}$ adjacent to $p_5$;
by hypothesis, $v$ is adjacent to $z$. Choose $a',b'\in L_{k-1}$ adjacent to $a,b$ respectively.
Consequently $a',b'$ are not adjacent to $z$, and so have no neighbours in $V(P)$; and in particular,
$a',b'$ are different from $v$ (although possibly $a'=b'$).
There is an even induced path between $v,a'$ with interior in $L_0\cup\cdots\cup L_{k-2}$,
and so the odd path $v\d p_5\d p_4\d p_3\d a\d a'$ is not induced, since its union with the previous path
would form an odd hole of length at least seven.
But $a'$ has no neighbour in $P$ (because $V(P)\subseteq A$), and $v$ is not adjacent to 
$a$ (because $G$ is triangle-free) and $v$ is not adjacent to $a'$ (because $L_{k-1}$ is stable), and it follows that $v$ is adjacent to $p_3$.
The same arguments applied to the path $v\d p_5\d p_1\d p_2\d b\d b'$ show that $v$ is adjacent to $p_2$;
yet not both of these are true since $G$ is triangle-free, a contradiction. This proves \ref{nbrz}.~\bbox

\bigskip

We deduce:

\begin{thm}\label{2balls}
If $H$ is a $2$-coverable graph and $z\in V(H)$ then $\chi(N^2_H(z))\le 5$.
\end{thm}
\Proof Since $H$ is $2$-coverable, there is a $1$-coverable graph $G$ and a 
stable levelling $(L_0\l L_k)$ in $G$ over $H$.
Let $A$ be the set of all vertices $v$ in $N^2_H(z)$ such that every neighbour of $v$ in $L_{k-1}$ is adjacent
to $z$, and let $B=N_H^2(z) \setminus A$. By \ref{nbrz}, 
$\chi(A)\le 2$, so we may assume (for a contradiction) that $\chi(B)> 3$.

Choose $z_0\in L_{k-1}$ adjacent to $z$. Since $N_G(z_0)$ is stable, it follows that
$\chi(B\setminus N_G(z_0))\ge 3$;
and so there is a $5$-hole $P$ with $V(P)\subseteq B$, such that $z_0$ has no neighbours in $P$.
Let $S_1\subseteq N_H(z)$
be minimal such that every vertex in $P$ has a neighbour in $S_1$. Each vertex in $P$ has a neighbour in $L_{k-1}$ nonadjacent to $z$,
and so there exists a minimal subset $S_2$ of $L_{k-1}\setminus N_G(z)$ such that every vertex in $P$ has a neighbour in $S_2$.
By \ref{pentcover}, $|S_1|=|S_2| = 3$. 
\\
\\
(1) {\em If $a_1\in S_1$ and $a_2\in S_2$ are joined by a three-edge path with interior in $V(P)$, then $a_1,a_2$
are adjacent. In particular, if $a_1\in S_1$ and $a_2\in S_2$ both have two neighbours in $V(P)$ and have a common neighbour in $V(P)$
then they have the same neighbours in $V(P)$.}
\\
\\
Let $a_1, a_2$ be adjacent to $p_1,p_2$ respectively, where $p_1p_2$ is an edge of $P$. 
If $a_1,a_2$ are not adjacent, then 
the path 
$z_0\d z\d a_1\d p_1\d p_2\d a_2$ is induced, and can be completed to an odd hole of length at least seven
via an even induced path between $z_0,a_2$ with interior in $L_0\cup\cdots\cup L_{k-2}$, which is impossible. This proves 
the first claim of (1). For the second, suppose that $a_1,a_2$ have a common neighbour in $V(P)$; then they are nonadjacent, and
so cannot be joined by a three-edge path with interior in $V(P)$, by the first claim. This proves (1).

\bigskip

Let $S_i=\{a_i,b_i,c_i\}$ for $i = 1,2$.
By \ref{pentcover}, for $i = 1,2$ we may assume that $a_i,b_i$ each have two neighbours in $V(P)$, 
and have no common neighbour in $V(P)$.
So one of $a_2,b_2$, say $a_2$, is adjacent to a neighbour of $a_1$ in $V(P)$, and hence 
$a_1,a_2$ have the same neighbours in $V(P)$, by the second claim of (1).
Therefore $b_2$ and $b_1$ have a common neighbour in $V(P)$, and so by the same argument, $b_1,b_2$
have the same neighbours in $V(P)$. If $c_1$ has two neighbours in $V(P)$, then it has a common neighbour in $V(P)$
with one of $a_2,b_2$, and so by the second claim of (1) it has the same neighbours in $V(P)$ as one of $a_2,b_2$, 
and hence the same as one
of $a_1,b_1$, which is impossible by the minimality of $S_1$. Thus $c_1$ has exactly one neighbour in $P$, and similarly
$c_2$ has exactly one neighbour in $P$, and the same neighbour as $c_1$. 

We may therefore assume that for $i = 1,2$, $a_i$ is adjacent to $p_2,p_4$ and $b_i$ to $p_3,p_5$, and $c_i$
to $p_1$. By the first claim of (1), it follows that $a_1$ is adjacent to $b_2,c_2$, and $b_1$ to $a_2,c_2$, and
$c_1$ to $a_2,b_2$.
But then the subgraph induced on
$$\{p_1,p_2,p_4,p_5, a_1,b_1,c_1,a_2,b_2,c_2,z\}$$
is isomorphic to the Gr\"{o}tzsch graph (with rim $a_2\d c_1\d p_1\d p_5\d b_1\d a_2$ and apex $a_1$), contradicting
\ref{myccover} since $G$ is $1$-coverable.
This proves \ref{2balls}.~\bbox

\bigskip

Now we complete the proof of \ref{mainthm}, which we restate:
\begin{thm}\label{mainthm3}
Every pentagonal graph is $58000$-colourable.
\end{thm}
\Proof
Define $n_1 = 581$, $n_2 = 10n_1-9$, and $n_3 = 10n_2-9$.
Suppose that there is a pentagonal graph $G_3$ with $\chi(G_3)\ge n_3$. By \ref{getlevelling}, there is a
stable levelling in $G_{3}$ over some graph $G_{2}$
with $\chi(G_{2})\ge n_{2}$.
Similarly there is a stable levelling in $G_2$ over some $G_1$ with $\chi(G_1)\ge n_1$. By \ref{2balls},
$\chi(N^2_{G_1}(v))\le 5$ for every vertex $v$
of $G_1$. By \ref{2rad} with $l=2$ and $\kappa_2=5$ it follows that $\chi(G_1)\le (40l+28)\kappa_2+40= 580$,
a contradiction.
Thus there is no such $G_3$, and hence every pentagonal graph has chromatic number at most $n_3-1 = 58000$.
This proves \ref{mainthm3}.~\bbox

\section{$7$-holes}
The previous result shows that triangle-free graphs with large chromatic number contain odd holes of length at 
least seven. But what if we ask for a hole of length {\em exactly} seven? We need to study this for an application
in a future paper.
Of course, the result is not true any more;
graphs with large chromatic number can have large girth. But we can still rescue something, by modifying the
foregoing proofs, with an additional hypothesis. We need to assume that in every induced subgraph with large chromatic
number, there is a vertex $v$ with $\chi(N^2(v))$ large. 

More precisely, let $\mathbb{N}$ denote the set of non-negative integers, and let 
$\phi:\mathbb{N}\rightarrow \mathbb{N}$ be a non-decreasing function.
We say a graph $G$ is {\em $\phi$-balled} if for every non-null induced subgraph $H$ of $G$,
there is a vertex $v\in V(H)$ such that $\chi(H)\le \phi(\chi(N^2_H(v)))$. We claim:

\begin{thm}\label{phiballs}
For every non-decreasing function $\phi:\mathbb{N}\rightarrow \mathbb{N}$, every 
triangle-free $\phi$-balled graph with no $7$-hole has chromatic number at most $\phi(\phi(\phi(2\phi(1)+1)))$.
\end{thm}
\Proof
As we said, the proof is a modification of the previous arguments (replacing the ``pentagonal'' condition in 
the definition of $n$-coverable and in \ref{pentcover}, \ref{nbrz} by a ``no $7$-hole'' condition),
so we just sketch it.
First, we modify the definition of a levelling; now we only consider levellings $(L_0, L_1, L_2)$ with at most 
three members. We define $n$-coverable as before (with this modification, and with ``no $7$-hole'' replacing ``pentagonal''). 
\ref{myccover} remains true under
this modification (to see this, observe that in the proof of \ref{myccover}, if $k=2$ then all the odd holes of length
at least seven found in that proof are in fact $7$-holes.) Also \ref{pentcover} still holds. 
To go further we need the following.
\\
\\
(1) {\em If $G$ is triangle-free and $\phi$-balled, and $\chi(G)>\phi(1)$, then $G$ has a $5$-hole.}
\\
\\
Choose a vertex $v$ such that $\chi(G)\le \phi(\chi(N^2(v)))$. Since $\chi(G)>\phi(1)$ and $\phi$ is non-decreasing,
it follows that $\chi(N^2(v))>1$ and in particular some two vertices in $N^2(v)$ are adjacent; and since $G$ is triangle-free
it follows that $G$ has a $5$-hole as required. This proves (1).

\bigskip

With the aid of (1), we can resuscitate 
a version of \ref{nbrz}, assuming that $G$ has no $7$-hole. At the start of the proof,
we have a set $A$ with $\chi(A)\ge 3$, and we deduce that there is a $5$-hole $P$ in $G[A]$. This is
no longer valid, because there might instead be a hole of length greater than seven. But if we replace $3$
by $\phi(1)+1$, then we can apply (1) to get a $5$-hole, and the remainder of the proof works.
The new version of \ref{nbrz} says that $\chi(A)\le \phi(1)$, with $A$ as before.
Similarly we can use (1) to obtain a form of \ref{2balls}; the new version says that if $H$ has no $7$-hole then
$\chi(N^2_H(z))\le 2\phi(1)+1$.

The proof of \ref{mainthm3} becomes easier, because we no longer need \ref{getlevelling} or \ref{2rad}. 
From \ref{2balls},
and the definition of $\phi$-balled, every 2-coverable triangle-free $\phi$-balled graph with no $7$-hole has chromatic number at most
$\phi(2\phi(1)+1)$. Consequently every $1$-coverable such graph has chromatic number at most 
$\phi(\phi(2\phi(1)+1))$ ; and so every triangle-free $\phi$-balled graph with no $7$-hole has chromatic number at most 
$\phi(\phi(\phi(2\phi(1)+1)))$. This proves \ref{phiballs}.~\bbox

\section{Long holes}

In this section we prove \ref{longoddhole} and \ref{longhole}. The first is implied by the next result with $m=2$:

\begin{thm}\label{longoddhole2}
Let $\ell\ge m\ge 2$ be integers, and let $G$ be a triangle-free graph with no odd hole of 
length at most $2m+1$ and no odd hole of length more than $2\ell+1$.
Then $\chi(G)< (3+4\ell)4^{\ell-m}-4\ell$.
\end{thm}
\Proof 
We proceed by induction on $\ell-m$. If $m=\ell$ then $G$ is bipartite and the result is true, so we assume that $m<\ell$.
Suppose that $\chi(G)\ge (3+4\ell)4^{\ell-m}-4\ell$. Then we may choose a levelling in $G$ with base of chromatic number at least 
$\chi(G)/2 \ge (6+8\ell)4^{\ell-m-1}-2\ell.$
Since $G$ has no odd cycle of length at most five, it follows that $N^2(v)$
is stable for every vertex $v$; and so by \ref{parentrule} with $\kappa=1$,
there is a stable levelling  $(L_0,L_1\l L_k)$ in $G$ with $\chi(L_k)\ge (3+4\ell)4^{\ell-m-1}-2\ell$,
and we may choose it such that $G[L_k]$ is connected.
It follows that $k\ge 3$.
For $0\le i\le k$ choose $s_i\in L_i$ such that $s_0\d s_1\c s_k$ is a path.
Since $\chi(L_k)>2l$ and $(L_k,s_{k-2}\d s_{k-1})$ is a lollipop,
\ref{lollipop} with $\kappa=1$ implies that there is a licking $(C,T_1)$ of this lollipop 
with 
$$\chi(C)\ge \chi(L_k)-2\ell= (3+4\ell)4^{\ell-m-1}-4\ell$$ 
and cleanliness at least $2\ell$. From the inductive hypothesis,
there is a $(2m+3)$-hole $P$ in $C$, with vertices $p_1\c p_{2m+3}\d p_1$ say. By \ref{licking} there is a licking $(P,T)$
of $(C,T_1)$. Let $T$ have vertices 
$$s_{k-2}\d s_{k-1}\d t_1\c t_r$$ 
say; thus $t_r$ has a neighbour in $P$, and since
the lollipop $(P,T)$ has cleanliness at least $2\ell$, it follows that $r\ge 2\ell$ and each of 
$s_{k-2}, s_{k-1}, t_1\l t_{2\ell-2}$ has distance at least three from $V(P)$. 

Now since $G$ has no odd cycle of length less than $2m+3$, it follows that every vertex of $G$ not in $P$
either has at most one neighbour in $P$, or has exactly two neighbours in $P$ with distance two in $P$. We may therefore
assume that $t_r$ is adjacent to $p_1$ and to no other vertex of $P$ except possibly $p_{2m+2}$.
For $i=3,4$, choose $a_i\in L_{k-1}$ adjacent to $p_i$. It follows that $a_3,a_4$ are nonadjacent
to $s_{k-2}, s_{k-1}, t_1\l t_{2l-2}$. Since $L_0\l L_{k-1}$ are stable, for $i = 3,4$ there is an even induced path
$R_i$ between $a_i$ and $s_{k-1}$ with interior in $L_0\cup\cdots\cup L_{k-2}$.
\\
\\
(1) {\em $a_4$ has a neighbour in $V(T)$.}
\\
\\
Because suppose not. Then $R_4\cup T$
is an induced path from $a_4$ to $t_r$, of length at least $r+2\ge 2l+2$.
But there is an odd induced path and an even induced path between $a_4$ and $t_r$ with interior in $V(P)$
(since $a_4$ has no neighbours in $P$ except $p_4$ and possibly $p_2,p_6$, and $t_r$ has no neighbours in 
$P$ except $p_1$ and possibly $p_{2m+2}$;  one of $a_4\d p_4\d p_3\d p_2\d p_1\d t_4$,
$a_4\d p_2\d p_1\d t_r$ is the desired odd path, and the even path goes the other way around $P$.) But then the
union of one of these paths with $R_4\cup Q$ is an odd hole of length at least $2l+4$, which is impossible. This proves (1).

\bigskip

Choose $i\le r$ minimum such that $t_i$ is adjacent to one of $a_4,a_3$. By (1), such a choice is possible. Since
$a_3, a_4$ are nonadjacent
to $s_{k-2}, s_{k-1}, t_1\l t_{2\ell-2}$, it follows that $i>2\ell-2$. Since $G$ has no odd cycle of length at most five,
$t_i$ is not adjacent to both $a_3,a_4$; let $t_i$ be adjacent to $a_h$ and not to $a_j$, where $\{h,j\}= \{3,4\}$.
Let $Q$ be a minimal path between $a_h, s_{k-1}$ with interior in $V(T)$. It follows that $Q$ has length at least $2\ell$. Consequently
$Q\cup R_h$ is a hole of length at least $2\ell+2$, and so it is even; and hence $Q$ is even. Now $a_h$ has 
no neighbour in $R_j$, since $a_h$ is not adjacent to the parent of $a_j$ (because $G$ has no $5$-holes)
and $a_h$ is nonadjacent to $s_{k-2}$ (because $(P,T)$ is a lollipop of cleanliness at least one). Thus
$$a_h\d p_h\d p_j\d a_j\d R_j\d s_{k-1}\d Q\d a_h$$
is an odd hole of length at least $2\ell+5$, which is impossible. This proves \ref{longoddhole2}.~\bbox 

Finally we turn to the proof of \ref{longhole}.
It follows from the next result.

\begin{thm}\label{longholelemma}
Let $\ell\ge 3$ and $\kappa\ge 1$ be integers, and let $G$ be a graph with no hole of length more than $\ell$, such that
$\chi(N(v)),\chi(N^2(v))\le \kappa$ for every vertex $v$.
Then $\chi(G)\le (2\ell-2)\kappa$.
\end{thm}
\Proof 
Suppose not; then there is a levelling $(L_0\l L_k)$ in $G$ with $\chi(L_k)> (\ell-1)\kappa$.
Let $C'$ be the vertex set of a component $C'$ of $G[L_k]$ with $\chi(C')> (\ell-1)\kappa$.
Since $\ell-1>1$, it follows that $k\ge 2$.
For $i = k-2,k-1$ choose $s_i\in L_i$, such that $s_{k-2}, s_{k-1}$ are adjacent and $s_{k-1}$ has a neighbour in $C'$.
Since $\chi(C')>(\ell-1)\kappa$ and $(V(C'),s_{k-2}\d s_{k-1})$ is a lollipop, by \ref{lollipop} there is a licking
$(C, T)$ of it with cleanliness at least $\ell-1$ and with $\chi(C)\ge \chi(C')-(\ell-1)\kappa>0$.
Choose $a\in L_{k-1}$ with a neighbour in $C$. Now $a$ might have neighbours in $T$, but since $(C, T)$
has cleanliness at least $\ell-1$, $a$ is nonadjacent to the first $\ell-1$ vertices of $T$. Let $P$ be an induced path between $s_{k-1}$
and $a$ with interior in $V(T)\cup C$; thus $P$ has length at least $l-1$. But $a,s_{k-1}$ are joined by an induced path
with interior in $L_0\cup \cdots\cup L_{k-2}$, and the union of this path with $P$ is a hole of length at least $\ell+1$,
a contradiction. This proves \ref{longholelemma}.~\bbox

\bigskip

We deduce \ref{longhole}, which we restate, slightly strengthened.
\begin{thm}\label{longhole2}
Let $\ell\ge 3$ be an integer, and let $G$ be a graph with no $5$-hole and no hole of length more than $\ell$. Then 
$$\chi(G)\le (2\ell-2)^{2^{\omega(G)-1}-1}.$$
\end{thm}
\Proof
We proceed by induction on $\omega(G)$. If $\omega(G)=1$ the result is true, so we assume $\omega(G)>1$.
Let 
$$n= (2\ell-2)^{2^{\omega(G)-2}-1}.$$
From the inductive hypothesis, every induced subgraph $H$ of $G$ with $\omega(H)<\omega(G)$ is $n$-colourable.
\\
\\
(1) {\em For every vertex $v$ of $G$, $\chi(N(v))\le n$, and $\chi(N^2(v))\le n^2$.}
\\
\\
The graph $G[N(v)]$ contains no clique of size $\omega(G)$,
and so is $n$-colourable. Let $A_1\l A_n$ be a partition of $N(v)$ into $n$ stable sets, and for $1\le i\le n$ let $B_i$
be the set of vertices in $N^2(v)$ with a neighbour in $A_i$. Suppose that there is a clique $C$ of cardinality $\omega(G)$
with $C\subseteq B_i$ for some $i$. Choose $a\in A_i$ with as many neighbours in $C$ as possible; then there exists $c'\in C$
nonadjacent to $a$, since $G$ has no $(\omega(G)+1)$-clique. Choose $a'\in A_i$ adjacent to $c'$; then from the
choice of $a$, there exists $c\in C$ adjacent to $a$ and not to $a'$. But then the subgraph induced on $\{v,a,a',c,c'\}$
is a $5$-hole, which is impossible. Thus there is no such clique $C$, and so $\chi(A_i)\le n$. Since this holds for
all $i$, it follows that $\chi(N^2(v))\le n^2$. This proves (1).

\bigskip
From (1) and \ref{longholelemma}, it follows that 
$$\chi(G)\le (2\ell-2)n^2 =  (2\ell-2) (2\ell-2)^{2^{\omega(G)-1}-2} = (2\ell-2)^{2^{\omega(G)-1}-1}.$$
This proves \ref{longhole2}.~\bbox

\end{document}